\documentclass[oneside]{article}
\usepackage{amsmath, amssymb, enumerate, graphicx, array}
\newtheorem{lemma}{Lemma}
\newtheorem{cor}{Corollary}
\newcommand{\tr}{\operatorname{tr\,}}
\title{The geometry of the set of real square roots of $\pm I_2$}
\author{{\bf Dr. V.N. Krishnachandran} \\ [0.2mm] 
Professor of Computer Applications \\
Vidya Academy of Science \& Technology \\
Thalakkottukara, Thrissur - 680501\\
email: {\tt krishnachandran.v.n@vidyaacademy.ac.in}}
\date{}
\allowdisplaybreaks
\begin{document}
\maketitle
\vspace{1.75cm}

\begin{abstract}

In this paper we study the geometry of the set of real square roots of $\pm I_2$. After some introductory remarks, we begin our study by deriving by quite elementary methods the forms of the real square roots of $\pm I_2$. We then discuss the interpretations of these square roots as transformations of the cartesian $(x,y)$-plane. To study the geometry of the set of square roots of $\pm I_2$ we consider a slightly more general  set of square matrices of order $2$ and show that these sets are hyperboloids of one sheet or hyperboloids of two sheets. From these general results we conclude that the set of involutory matrices of order 2 is a hyperboloid of one sheet and the set of skew-involutory matrices of order 2 is a hyperboloid of two sheets. The relations between the geometrical properties of the hyperboloids and the set of square roots of $I_2$ are also investigated. We then proceed to obtain the forms of the involutory matrices of order 2 by more advanced methods. We have considered two approaches: in the first approach we use the concept of a function of a matrix and in the second approach we use concepts of split-quaternions.
\end{abstract}
\newpage
\ 
\vspace{1in}
\tableofcontents
\newpage
\section{Introduction}
\subsection{Matrix square roots}
Let $P$ be a square matrix over a field $F$. A square matrix $R$ over $F$ such that $R^2=P$ is called a {\em square root of $P$}.

The term ``matrix" was coined by James Joseph Sylvester in 1850 and matrix algebra was developed by Arthur Cayley in his ``Memoir on the Theory of Matrices'' published in 1858 \cite{Cayley}. In this memoir, Cayley considered matrix square roots also. 

An arbitrary  square matrix may have no, or a finite number of, or an infinite number of square roots.  For example,  it can be seen that:
\begin{itemize}
\item
$\begin{bmatrix} 0 & 1 \\ 0 & 0\end{bmatrix}$ has {\em no} square root;
\item
$\begin{bmatrix} 1 & 0 \\ 0 & 0 \end{bmatrix}$ has exactly {\em two} square roots, namely, $\begin{bmatrix} \pm 1 &  \phantom{\pm} 0 \\ \phantom{\pm} 0 &  \phantom{\pm} 0 \end{bmatrix}$;
\item
$\begin{bmatrix} 1 & 0 \\  0 & 4 \end{bmatrix}$ has exactly {\em four} square roots, namely, $\begin{bmatrix} \pm 1 &  \phantom{\pm} 0 \\  \phantom{\pm} 0 & \pm 2\end{bmatrix}$;
\item
$\begin{bmatrix} 4 & 0 \\ 0 & 4 \end{bmatrix}$ has an {\em infinite number} of square roots some of which are given by $\begin{bmatrix}2 \cos \phi  & \phantom{-}2\sin\phi  \\2 \sin\phi  & - 2\cos\phi \end{bmatrix}$ for arbitrary values of $\phi$. 
\end{itemize} 

The matrix square root is one of the most  commonly occurring matrix functions,  and it arises in several contexts: for example in the 
matrix sign function, the generalized eigenvalue problem, the polar decomposition, and the geometric mean problem. There are also a 
variety of methods for computing the matrix square root (see, for example, \cite{Nigham}, \cite{Jean}). When the order of the matrix is $2$, the matrix becomes particularly simple and there are special formulas and techniques to compute the square roots of matrices of order $2$ (see, for example, \cite{Sullivan},\cite{Sam}).  
\subsection{Involutory and skew-involutory matrices}
In this paper, we study the square roots of a very special and  a very interesting square matrix of order $2$, namely, the identity matrix $I_2= \begin{bmatrix} 1 & 0 \\ 0 & 1\end{bmatrix}$ of order $2$. In the literature, square roots of the identity matrix $I$ is called an {\em involutory matrix} and the square roots of $-I$ are called {\em skew-involutory matrices}. Thus this paper is concerened with the geometry of the set of involutory and skew-involutory matrices of order 2. 

In Section \ref{S2}, we first find {\em all} the square roots of  the identity matrix $I_2$  by elementary methods. In Section \ref{S3} we discuss the geometrical interpretations of these square roots as transformations of the cartesian $xy$-plane. 
\section{The square roots of $I_2$: Elementary approach}\label{S2}
Let $R=\begin{bmatrix} a & b \\ c& d \end{bmatrix}$ be a square root of $I_2$. The defining condition $R^2 = I_2$ now implies the following relations:
\begin{align}
a^2 + bc & = 1 \label{Eq1}\\
(a+d)b & = 0 \label{Eq2}\\
(a+d)c & =0 \label{Eq3}\\
d^2 + bc & =1 \label{Eq4}
\end{align}
\subsubsection*{Case 1}
If $a+d\ne 0$ then from Eq.\eqref{Eq2} and  Eq.\eqref{Eq3} we get $b=0$ and $c=0$. Also from Eq.\eqref{Eq1} and Eq.\eqref{Eq4} we have $a=\pm 1$ and $d=\pm 1$. The condition $a+d\ne 0$ will be stisfied only when $a=d=1$ or $a=d=-1$. Thus, if $a+d\ne 0$, we have: 
\begin{enumerate}[{\hspace{\parindent}\bf {1}(a)}]
\item
$R=\begin{bmatrix}1 & 0 \\ 0 & 1 \end{bmatrix}=I_2$ 
\item
$ R=\begin{bmatrix}-1 & \phantom{-}0 \\ \phantom{-}0 & -1 \end{bmatrix}=-I_2$.
\end{enumerate}
\subsubsection*{Case 2}
Now assume that $a+d=0$ so that $d=-a$. Then Eq.\eqref{Eq2} and Eq.\eqref{Eq3} are obviously satisfied and the other two equations become equivalent to each other. We have to distinguish two cases:

\subsubsection*{Case 2(i)}
Let $a^2=1$ so that $bc=0$ and so $b=0$ or $c=0$. This case yields the following square roots of $I_2$:
\begin{enumerate}[{\bf {\hspace{\parindent} 2(i)}(a) }]
\item
$R= \begin{bmatrix}\phantom{-}1 & \phantom{-}b \\ \phantom{-}0 & -1 \end{bmatrix}$ (See Figure \ref{Fig1})
\item
$R=\begin{bmatrix}-1 & \phantom{-}b \\ \phantom{-}0 & \phantom{-}1 \end{bmatrix}$ (See Figure \ref{Fig2})
\item
$R=\begin{bmatrix}\phantom{-}1 & \phantom{-}0 \\\phantom{-} c & -1 \end{bmatrix}$ (See Figure \ref{Fig3})
\item
$R=\begin{bmatrix}-1 & \phantom{-}0 \\ \phantom{-}c & \phantom{-}1 \end{bmatrix}$ (See Figure \ref{Fig4})
\end{enumerate}
\subsubsection*{Case 2(ii)}
Let $a^2 \ne 1$ so that $b\ne 0$ and $c\ne 0$. Now from Eq.\eqref{Eq1} we have
$$
c=\frac{1-a^2}{b}.
$$
This yields the following general expression for square roots of $I_2$:
\begin{enumerate}[{\bf \hspace{\parindent} {2(ii)}(a)}]
\item
$R=\begin{bmatrix} a & \phantom{-}b \\ \dfrac{1-a^2}{b} & -a \end{bmatrix}.$
\end{enumerate}

This shows that the set of square roots of $I_2$ is a two-parameter family of matrices.
\section{Geometrical interpretations of the square roots of $I_2$}\label{S3}
\subsection{Matrices of order $2$ as transformations of a plane}
The two-dimensional vector space $\mathbb R^2$ over $\mathbb R$ can be geometrically realised as the cartesian $xy$-plane. If $P(x, y)$ is any point in this plane, we may represent it by the column vector ${\mathbf x} = \begin{bmatrix} x \\ y \end{bmatrix}$. 

Let $T=\begin{bmatrix} a & b \\ c & d \end{bmatrix}$ be any square matrix of order $2$ over $\mathbb R$. Then the mapping $T:\mathbb R^2 \longrightarrow \mathbb R^2$ defined by $\mathbf x\mapsto  T\mathbf x$ defines the following   transformation of the $xy$-plane:
\begin{equation}\label{Trans}
(x,y)\mapsto (ax+by, cx+dy).
\end{equation}

Let $S$ be another square matrix of order $2$, then the transformation defined by the product $ST$ can be expressed as
$$
\mathbf x\mapsto  T\mathbf x \mapsto S(T \mathbf x).
$$
This represents successive applications of two transformations: first, the transformation represented by $T$ and then the transformation represented by $S$. 

A transformation such as one given by Eq.\eqref{Trans} can be described in terms of elementary geometric transformations like, rotation, magnification, translation, point and line reflections, etc. In fact if $T=I_2$, the map $\mathbf x\mapsto  T\mathbf x$ is the identity map of  $\mathbb R^2$ onto itself. 
\subsection{Square roots of $I_2$ as transformations of a plane}
In this section we describe the various square roots of $I_2$ in terms of these elementary geometric transformations.  As aids to visualisation we have also included diagrams illustrating the transformations. In the diagrams, the maps from $P$ to $Q$ and from $Q$ to $P$ are both given by $T$ so that $T^2$ is the map from $P$ to $P$ which is the identity mapping.
\subsubsection*{Case 1}
\begin{figure}[!h]
\begin{center}
\includegraphics{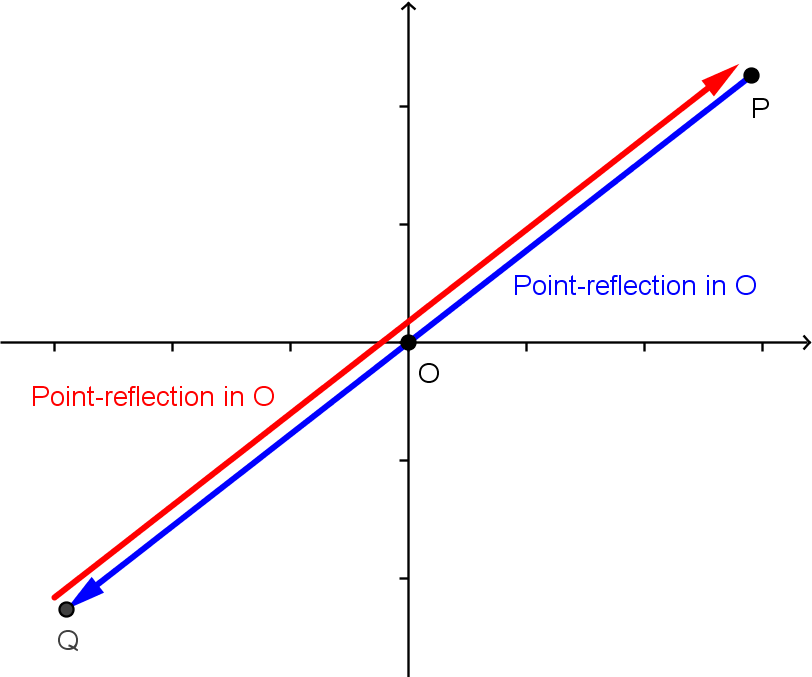}
\caption{Geometrical interpretation the square root $-I_2$ of $I_2$ (Case 1(b))}\label{Fig8}
\end{center}
\end{figure}

\begin{enumerate}[{\bf {\hspace{\parindent} 1}(a) }]
\item
As already mentioned, the square root $I_2$ of $I_2$ is the identity mapping of $\mathbb R^2$.
\item
The matrix $-I_2$ represents the point-reflection about the coordinate origin  (see Figure \ref{Fig8}).
\end{enumerate}
\subsubsection*{Case 2(i)(a)}
\begin{figure}
\begin{center}
\includegraphics{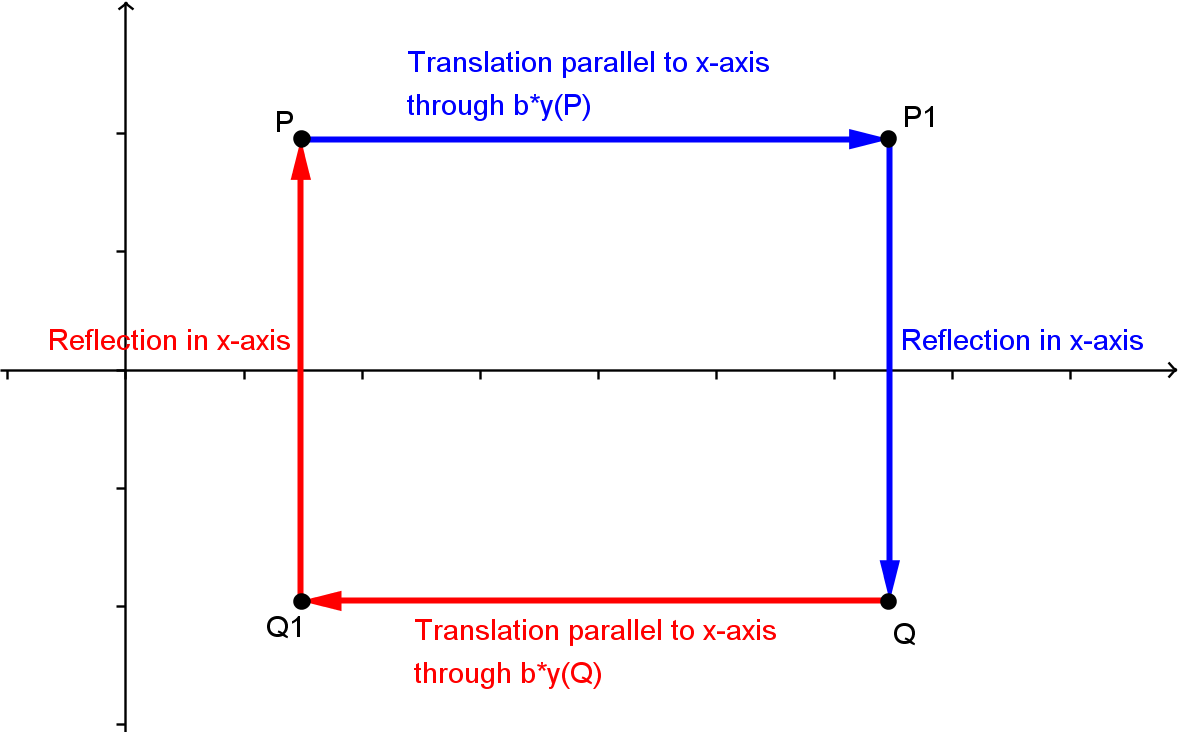}
\caption{Geometrical interpretation the square root  of $I_2$ (Case 2(i)a) }\label{Fig1}
\end{center}
\end{figure}

Next we have
$$
\begin{bmatrix}\phantom{-}1 & \phantom{-}b \\ \phantom{-}0 & -1 \end{bmatrix} = 
\begin{bmatrix}\phantom{-}1 & \phantom{-}0 \\ \phantom{-}0 & -1 \end{bmatrix}
\begin{bmatrix}\phantom{-}1 & \phantom{-}b \\ \phantom{-}0 & \phantom{-1}1 \end{bmatrix}
$$
The latter factor represents the map 
$$
(x,y)\mapsto (x+by, y)
$$
which is a translation of the point $P(x,y)$ parallel to $x$-axis through a distance of $by$. The former factor represents the map
$$
(x,y)\mapsto (x, -y)
$$
which is a refection  in the $x$-axis. The transfromations, and their repeated applications, are shown in Figure~\ref{Fig1}.
\subsubsection*{Case 2(i)(b)}
\begin{figure}[!t]
\begin{center}
\includegraphics[width=12cm]{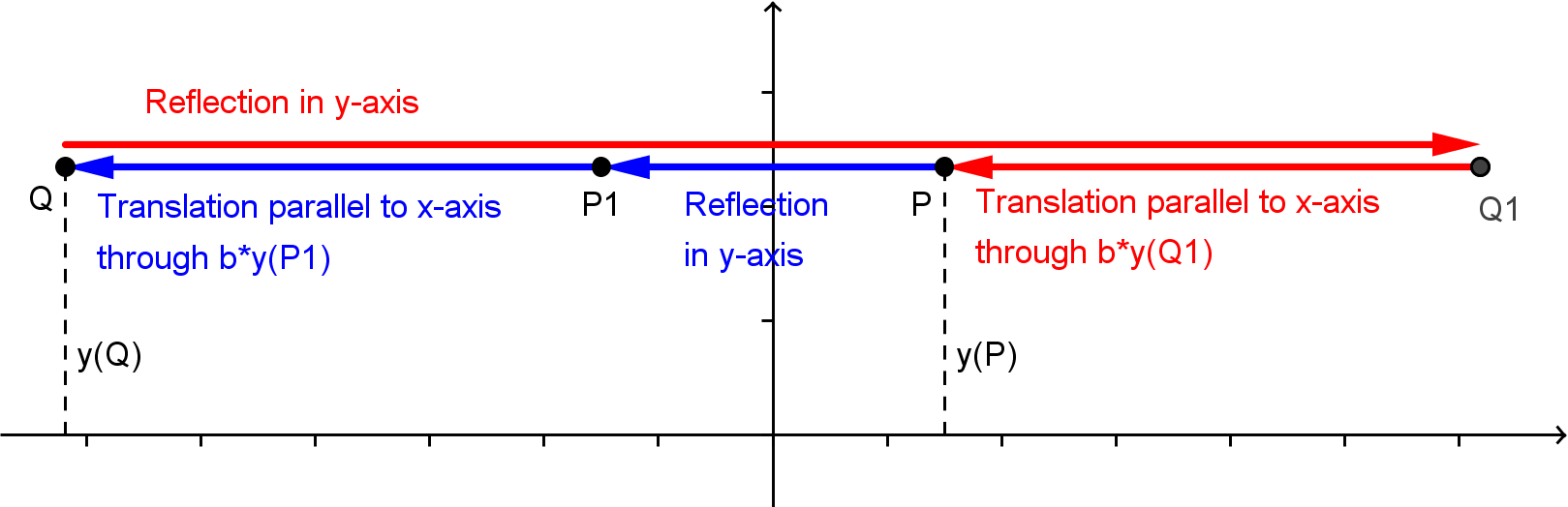}
\caption{Geometrical interpretation the square root  of $I_2$ (Case 2(i)b)}\label{Fig2}
\end{center}
\end{figure}

Next we have
$$
\begin{bmatrix} -1 & \phantom{-}b \\ \phantom{-}0 &  \phantom{-}1 \end{bmatrix} = 
\begin{bmatrix}\phantom{-}1 & \phantom{-}b \\ \phantom{-}0 & \phantom{-}1 \end{bmatrix}
\begin{bmatrix}-1 & \phantom{-}0 \\ \phantom{-}0 &  \phantom{-}1 \end{bmatrix}
$$
The latter factor represents the map 
$$
(x,y)\mapsto (-x, y)
$$
which is a refection  in the $y$-axis. The former factor represents the map
$$
(x,y)\mapsto (x+by,y)
$$
which is a translation of the point $P(x,y)$ parallel to $x$-axis through a distance of $by$.
 The transfromations, and their repeated applications, are shown in Figure~\ref{Fig2}. 
\begin{figure}[!h]
\begin{center}
\includegraphics{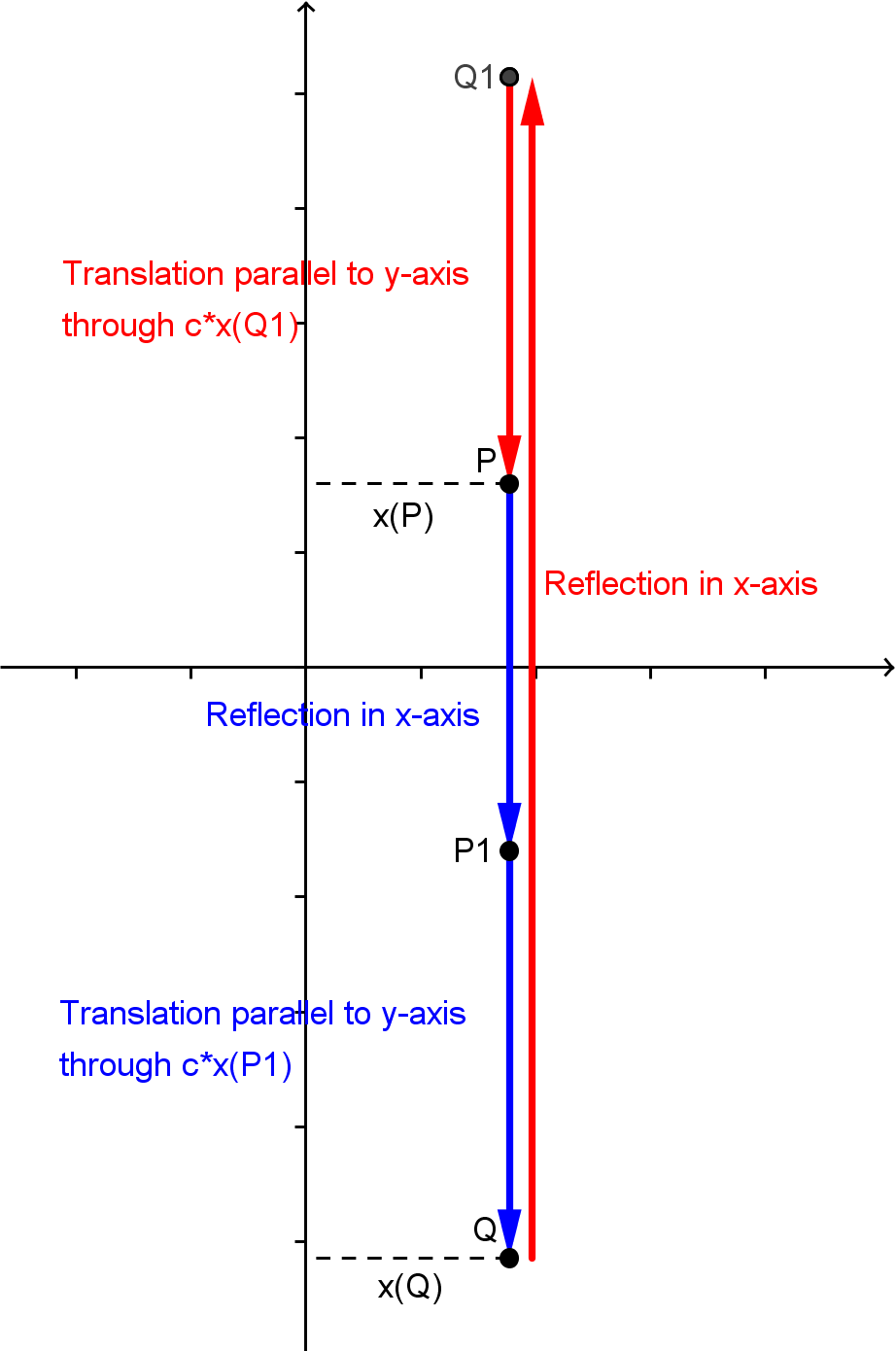}
\caption{Geometrical interpretation the square root  of $I_2$ (Case 2(i)c) }\label{Fig4}
\end{center}
\end{figure}
\begin{figure}[!h]
\begin{center}
\includegraphics{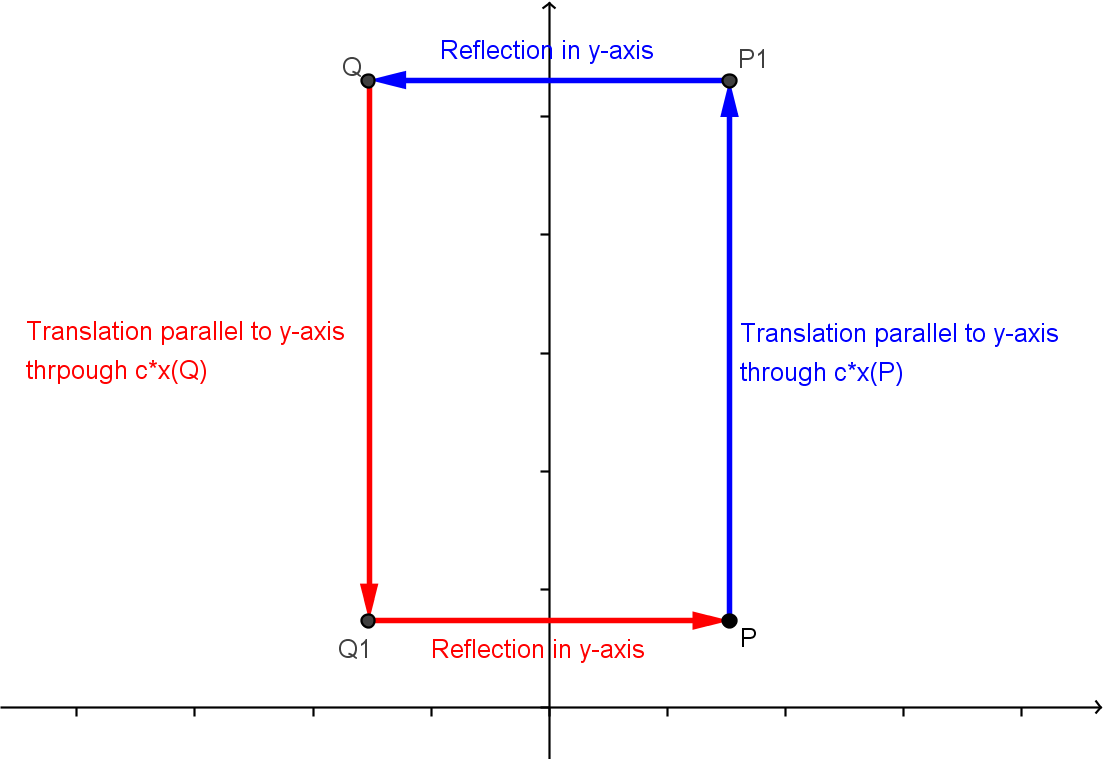}
\caption{Geometrical interpretation the square root  of $I_2$ (Case 2(i)d) }\label{Fig3}
\end{center}
\end{figure}

\subsubsection*{Case 2(i)(c)}
Next we have
$$
\begin{bmatrix} \phantom{-}1 & \phantom{-}0 \\ \phantom{-}c &  -1 \end{bmatrix} = 
\begin{bmatrix}\phantom{-}1 & \phantom{-}0 \\ \phantom{-}c & \phantom{-}1 \end{bmatrix}
\begin{bmatrix}\phantom{-}1 & \phantom{-}0 \\ \phantom{-}0 &  -1 \end{bmatrix}
$$
The latter factor represents the map 
$$
(x,y)\mapsto (x, -y)
$$
which is a refection  in the $x$-axis. The former factor represents the map
$$
(x,y)\mapsto (x,cx+y)
$$
which is a translation of the point $P(x,y)$ parallel to $y$-axis through a distance of $cx$.
 The transfromations, and their repeated applications, are shown in Figure~\ref{Fig4}. 
\subsubsection*{Case 2(i)(d)}
Next we have
$$
\begin{bmatrix} -1 & \phantom{-}0 \\ \phantom{-}c &  \phantom{-}1 \end{bmatrix} = 
\begin{bmatrix} -1  & \phantom{-}0 \\ \phantom{-}0 & \phantom{-}1 \end{bmatrix}
\begin{bmatrix}\phantom{-}1 & \phantom{-}0 \\ \phantom{-}c &  \phantom{-}1 \end{bmatrix}
$$
The latter factor represents the map 
$$
(x,y)\mapsto (x, cx+y)
$$
which is a translation of the point $P(x,y)$ parallel to $y$-axis through a distance of $cx$.
The former matrix represents 
$$
(x,y)\mapsto (-x,y)
$$
which is a refection  in the $y$-axis. The transfromations, and their repeated applications, are shown in Figure~\ref{Fig3}. 
\subsubsection*{Case 2(ii): The general case}
\begin{figure}[t]
\begin{center}
\includegraphics[height=12cm]{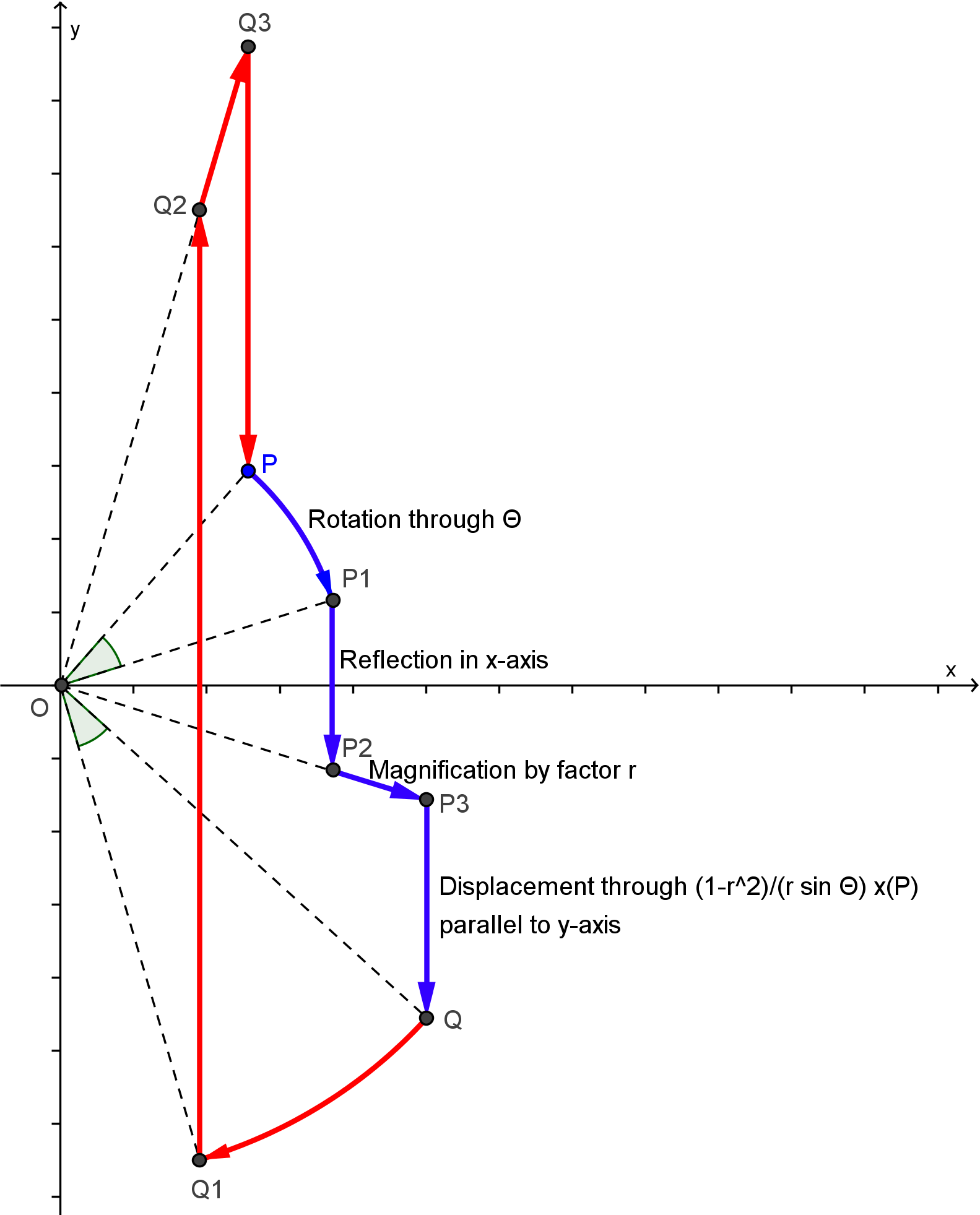}
\caption{Geometrical interpretation of the general square root of $I_2$}\label{Fig6}
\end{center}
\end{figure}

Now let us consider the general square root of $I_2$ given by 
$$
X=\begin{bmatrix} a & b \\ \frac{1-a^2}{b} & -a \end{bmatrix}.
$$
Letting $a=\rho\cos \phi$ and $b=\rho \sin\phi$, the matrix $X$ can be written as
\begin{equation}\label{R1R2}
R=R_3R_2R_1 + R_4
\end{equation}
where
$$
R_1=\begin{bmatrix}\phantom{-}\cos \phi & \sin\phi  \\ - \sin\phi  & \cos \phi  \end{bmatrix}, \quad
R_2 = \begin{bmatrix} 1 & \phantom{-}0 \\ 0 & -1 \end{bmatrix},\quad
R_3 = \rho I_2, \quad
R_4 = \begin{bmatrix} 0 & 0 \\ \frac{1-\rho^2}{\rho \sin\phi} & 0 \end{bmatrix}.
$$
The matrices $R_1, \ldots, R_4$ have the following geometrical interpretations (see Figure \ref{Fig6}). 
\begin{itemize}
\item
$R_1$ represents a rotation through the angle $\phi  $ in the clockwise direction.
\item
$R_2$ represents a reflection in the $x$-axis.
\item
 $R_3$  represents a magnification by a factor $\rho$.
\item
We have 
$$
R_4\begin{bmatrix}x \\ y\end{bmatrix} = \begin{bmatrix} 0  \\ 
\frac{1-\rho^2}{\rho\sin\phi} x \end{bmatrix}.
$$
Hence the term ``$+ \,R_4$'', that is, the addition of $R_4$ to $R_1R_2R_3$, in Eq.\eqref{R1R2} represents
  a translation through a distance $\frac{1-\rho^2}{r\sin\phi} x $ parallel to $y$-axis.
\end{itemize}
\section{A special case of the general case}
As a special case, let us consider those square roots of $I_2$ which are symmetric. The square roots discussed in Case 1 are clearly symmetric and those  discussed in Case 2 are symmetric only if $b=0$ or $c=0$ in which case they reduce to the matrices given in Case 1. 

The square roots discussed in the general case (that is, Case 3) are symmetric if 
$$
\frac{a^2-1}{b} =b;
$$
that is 
$$
a^2+b^2=1.
$$
Hence in  such a case we can 
take $a=\cos \phi $ and $b=\sin\phi $ for some $\phi $. Then the square root of $I_2$ takes the form
\begin{equation}\label{Householder}
X=\begin{bmatrix}\cos \phi & \phantom{-}\sin \phi \\
\sin\phi  & -\cos\phi  \end{bmatrix}
\end{equation}

Choosing a Pythagorean triple $(r,s,t)$ and letting $\cos \phi =\frac{r}{t}$ and $\sin\phi = \frac{s}{t}$, we have
$$
X=\frac{1}{t}\begin{bmatrix} r & \phantom{-}s \\ s & -r \end{bmatrix}.
$$
This case has been discussed in \cite{Gazette}.

\subsection{Geometrical interpretation: Rotation followed by reflection}
\begin{figure}[!h]
\begin{center}
\includegraphics{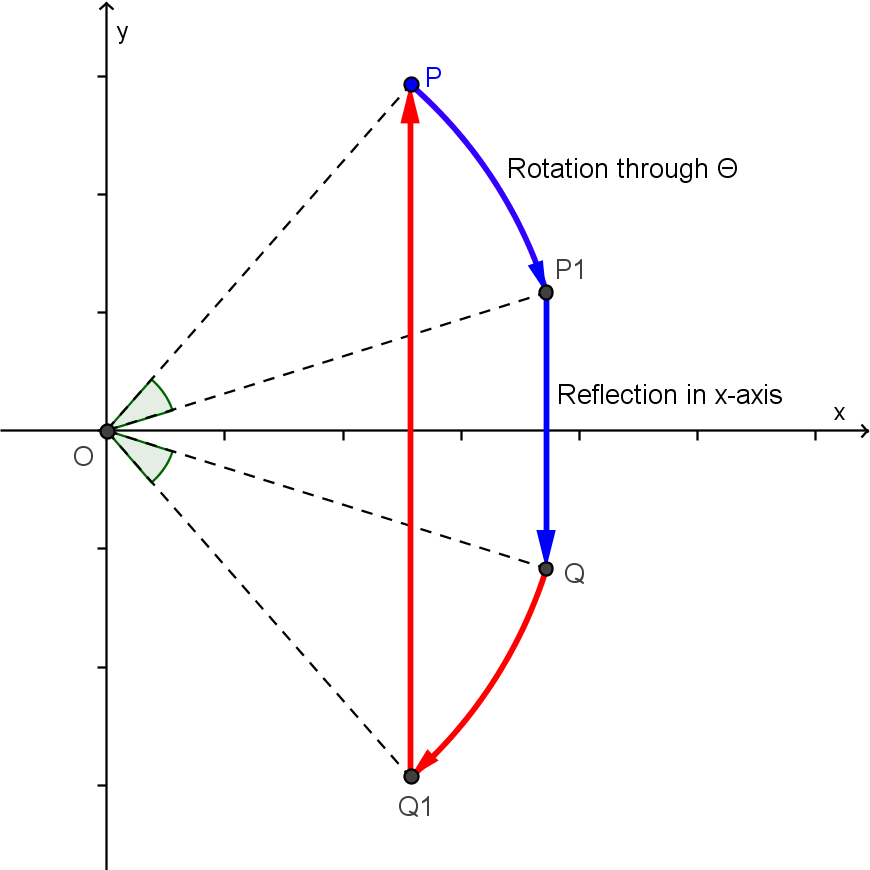}
\caption{Geometrical interpretation of a special square root of $I_2$}\label{Fig5}
\end{center}
\end{figure}

This square root has an interesting and obvious geometrical interpretation. 
The matrix can be written as
$$
R=R_2R_1
$$
where
$$
R_2= \begin{bmatrix}1 & \phantom{-}0  \\
0 & -1 \end{bmatrix}, \quad R_1 = \begin{bmatrix}\phantom{-}\cos \phi & \sin \phi \\
- \sin\phi & \cos\phi \end{bmatrix}.
$$
Here $R_2$ represents the reflection in the $x$-axis of points in the $xy$-cartesian plane and $R_1$ represents a rotation of points clockwise through an angle $\phi$ about the origin of the coordinate system (see Figure \ref{Fig5}).  Thus this special square root is a combination of rotation and reflection.
\subsection{Geometrical interpretation: Householder transformation}
The matrix given by Eq.\eqref{Householder} is a Householder matrix. 
\begin{figure}[!t]
\begin{center}
\includegraphics{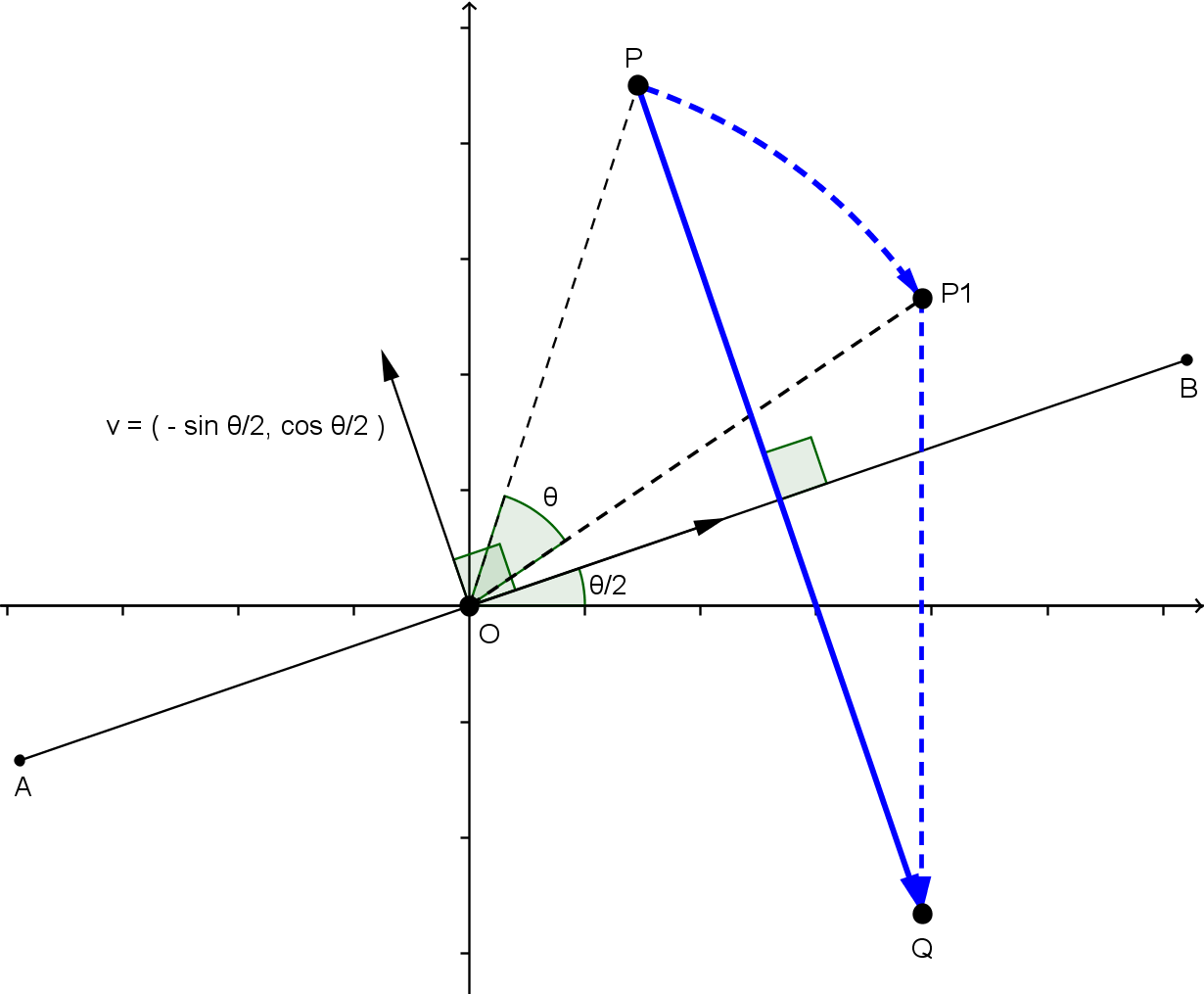}
\caption{A special square root of $I_2$ as a Householder transformation}\label{SquareRootOfIdentityMatrix10}
\end{center}
\end{figure}
Let us recall the definition of a Householder matrix. Let $\mathbf v$ be a column unit vector in $\mathbb C^n$ or $\mathbb R^n$ and $\mathbf v^H$ the Hermitian transpose of $\mathbf v$. Let $I$ be the identity matrix of order $n$. Then a matrix of the form
$$
P = I - 2 \mathbf v \mathbf v^H
$$
is called a Householder matrix. It is interesting to note that the identity matrix is not a Householder matrix. The map $\mathbf x \mapsto P\mathbf x$ is called a Householder transformation (see \cite{Householder} and \S 13.3 \cite{Datta}). Householder matrices are extensively used in numerical linear algebra. 

The matrix $P$ arises in the following way: Let $\mathbf x$ represent a point in space. Then the vector $P\mathbf x$ is the reflection of the point $\mathbf x$ in the hyperplane orthogonal to the unit vector $\mathbf v$.

Now to see that $X$ given by Eq.\eqref{Householder} is a Householder matrix, 
let
$$
\mathbf v = \begin{bmatrix} -\sin \frac{\phi}{2} \\ \phantom{-}\cos\frac{\phi}{2}\end{bmatrix}
$$
Then
$$
\mathbf v^H = \begin{bmatrix} -\sin\frac{\phi}{2} & \cos\frac{\phi}{2}\end{bmatrix}
$$
Now, we have
\begin{align*}
P & = I - 2 \mathbf v\mathbf v^H \\
& = \begin{bmatrix} 1 & 0 \\ 0 & 1 \end{bmatrix} - 2 \begin{bmatrix}\sin^2\frac{\phi}{2} & -\sin\frac{\phi}{2}\cos\frac{\phi}{2}  \\
\sin\frac{\phi}{2}\cos\frac{\phi}{2} & \cos^2\frac{\phi}{2}
\end{bmatrix}  \\
& = \begin{bmatrix}\cos \phi & \phantom{-}\sin \phi \\ \sin \phi & - \cos \phi \end{bmatrix} \\
& = X \qquad \text{(see Eq.\eqref{Householder})}.
\end{align*}
\section{Square roots of $-I_2$}
Since $I_2$ and $-I_2$ are closely related, in this section we briefly discuss the square roots of $-I_2$. 

As in Section 2, let $R=\begin{bmatrix} a & b \\ c& d \end{bmatrix}$ be a square root of $-I_2$ so that $R^2=-I_2$. So we must have 
\begin{align}
a^2+bc & = -1 \label{N1}\\
(a+d)b & =\phantom{-1}0 \label{N2}\\
(a+d)c & = \phantom{-1} 0 \label{N3}\\
d^2 + bc & = -1 \label{N4}
\end{align}
Proceeding as in Section 2, it can be seen that the real square roots of $-I_2$ are of the form
$$
R=\begin{bmatrix} a & \phantom{-}b \\ -\left(\frac{1+a^2}{b}\right) & -a 
\end{bmatrix}
$$
where $a$ and $b\ne 0$ are real numbers. The matrix $R$ as given above is a square root of $-I_2$ even in the case when $a$ and $b$ complex numbers. If we allow $a$ to be a complex number, then we must have $a^2+1\ne 0$.

The following are some special complex square roots of $-I_2$. All of them can be obtained as in Section 2.
$$
iI_2, \quad -iI_2,\quad
\begin{bmatrix} \phantom{-}i & \phantom{-}b \\ \phantom{-}0 & -i \end{bmatrix},\quad
\begin{bmatrix} -i & \phantom{-}b \\ \phantom{-}0 & \phantom{-}i \end{bmatrix},\quad
\begin{bmatrix} \phantom{-}i & \phantom{-}0 \\ \phantom{-}c & -i \end{bmatrix},\quad
\begin{bmatrix} -i & \phantom{-}0 \\ \phantom{-}c & \phantom{-}i \end{bmatrix}.
$$

%
%
\section{$S(\alpha, \beta)$: A generalisation of the set of square roots of $I_2$}\label{SectionA}
The set of all square roots of $I_2$ is  the set of all matrices $X$ satisfying the condition $X^2=I_2$. Similar such sets have been considered in the literature. 

For example, the geometry of the set $E_1$  of idempotent matrices of order $2$ has been  studied briefly in \cite{Krish01} and extensively in \cite{Krish02}. (Recall that a matrix $X$ is an idempotent if $X^2=X$.)  The set $E_1$ has been shown to be a hyperboloid of one sheet embedded in the four-dimensional space $\mathbb R^4$. The set of nilpotent matrices, that is, the set of matrices $X$ for which $X^2=0$,  is known to be a right circular cone. Some of these results on the set of idempotent matrices of order $2$ have been extended to sets of arbitrary order idempotent matrices (see \cite{Krish03}). 

The conditions defining the various sets of matrices described above are all 
special cases of  the following general condition: 
\begin{equation}\label{Eq5}
X^2 - \alpha X + \beta I_2=0.
\end{equation}
We shall now consider the geometry of the set of matrices $X$ satisfying Eq.\eqref{Eq5}. There are two special matrices $X$ satisfying Eq.\eqref{Eq5}. 
Let $t$ be a solution of the equation
\begin{equation}\label{Eq6}
t^2 - \alpha t +\beta =0.
\end{equation}
 then it can be seen that the matrix  $X=tI_2$ satisfies Eq.\eqref{Eq5}. The two solutions of Eq.\eqref{Eq6} yield two such matrices. In the sequel, when we consider solutions of Eq.\eqref{Eq5} we exclude these special solutions. 

To simplify the discussions, we introduce the following notation:
\begin{equation}\label{Eq7}
S(\alpha, \beta)= \left\{ X : X^2 -\alpha X + \beta I_2 =0 \text{ and } X\ne tI_2 \text{ for some } t\in \mathbb R\right\}.
\end{equation}
\subsection{Special cases}
We have the following interesting special cases:
\begin{itemize}
\item
$S(0,-1)$: This is the set of square roots of $I_2$ excluding the two roots $\pm I_2$.
\item
$S(1,0)$: This is the set $E_1$ of idempotent matrices of order 2 excluding $I_2$ and the zero matrix, studied in \cite{Krish02}.
\item
$S(0,0)$: This is the set of nilpotent matrices of order 2 excluding the zero matrix. 
\end{itemize}
\section{A simple characterization of the set $S(\alpha, \beta)$}
We now let 
$$
X=\begin{bmatrix} x_1 & x_2 \\ x_3 & x_4 \end{bmatrix}
$$
and identify $X$ with the point $(x_1, x_2, x_3, x_4)$ in the four-dimensional space $\mathbb R^4$.  
\begin{lemma}\label{Lemma1}
Let $X$ and $S(\alpha, \beta)$ be as above. Then $X\in S(\alpha, \beta)$ if and only if the following two conditions are satisfied:
\begin{align}
x_1+x_4 & = \alpha \label{Eq6a}\\
x_1x_4 - x_2x_3 & =\beta\label{Eq6b}
\end{align}
\end{lemma}
\paragraph{Proof.} 
Let $X$ be in $S(\alpha, \beta)$ and so $X$ satisfies Eq.\eqref{Eq5}. 

If we carry out the indicated operations in Eq.\eqref{Eq5} and equate the corresponding elements, we get the following system of equations:
\begin{align}
x_1^2+x_2x_3 - \alpha x_1 + \beta & =0 \label{Eq7}\\
(x_1+x_4)x_2 -\alpha x_2 & = 0 \label{Eq8}\\
(x_1+x_4)x_3 - \alpha x_3 & = 0 \label{Eq9}\\
x_4^2 +x_2x_3 -\alpha x_4 +\beta & =0 \label{Eq10}
\end{align}
Assume that $x_1+x_4\ne 2\alpha$. 
Subtracting Eq.\eqref{Eq10} from Eq.\eqref{Eq7}, we have
$$
(x_1-x_4)(x_1+x_4-\alpha)=0
$$
which implies that $x_1-x_4=0$. Also, from the remaining two equations we get $x_2=0$ and $x_3=0$. All these imply that $X=tI_2$ where $t=x_1=x_4$ which is a contradiction. Thus we must have $x_1+x_4=\alpha$. 
Further, replacing $x_1-\alpha$ by $- x_4$ in Eq.\eqref{Eq7} and simplifying we get
$$
x_1x_4-x_2x_3=\beta.
$$
Thus if $X\in S(\alpha, \beta)$ then Eq.\eqref{Eq6a} and Eq.\eqref{Eq6b} are satisfied. 

Conversely, it is straight forward to verify that if for a matrix $X$ the equations Eq.\eqref{Eq6a} and Eq.\eqref{Eq6b} are satisfied, then the Eqs.\eqref{Eq7}--\eqref{Eq10} are also satisfied so that $X$ is in $S(\alpha, \beta)$. $\blacksquare$
%

Since, $x_1+x_4=\tr (X)$, the trace of $X$, and $x_1x_4-x_2x_3=\det (X)$, the determinant of $X$, Lemma \ref{Lemma1} implies that $X\in S(\alpha,\beta)$ if and only if
$$
\tr (X)=\alpha, \quad \det(X)=\beta.
$$
Thus we have
$$
S(\alpha, \beta) = \{ X : \tr(X)=\alpha\} \cap \{ X : \det(X)=\beta \}
$$
To simplify the discussions we introduce the following notation:
$$
P(\alpha) = \{ X : \tr(X)=\alpha\}.
$$
\section{The hyperplane $P(\alpha)$ in $\mathbb R^4$}
As we have already indicated in the beginning of Section \ref{SectionA}, the square matrix $X=\begin{bmatrix} x_1 & x_2 \\ x_3 & x_4 \end{bmatrix}$ of order $2$ can be identified with the point $(x_1, x_2, x_3, x_4)$ in the four-dimensional Euclidean space $\mathbb R^4$. This correspondence is a bijection. Now, the set of points $X$ in $\mathbb R^4$ which satisfies the equation $\tr(X)=\alpha$  defines a hyperplane in this four-dimensional space. A hyperplane in $\mathbb R^4$ is a three-dimensional Euclidean space. Thus the set of matrices $S(\alpha, \beta)$ is actually a set of points in a three-dimensional Euclidean space $P(\alpha)$. We now investigate the geometry of this set of points. 
\begin{figure}[!h]
\begin{center}
\includegraphics{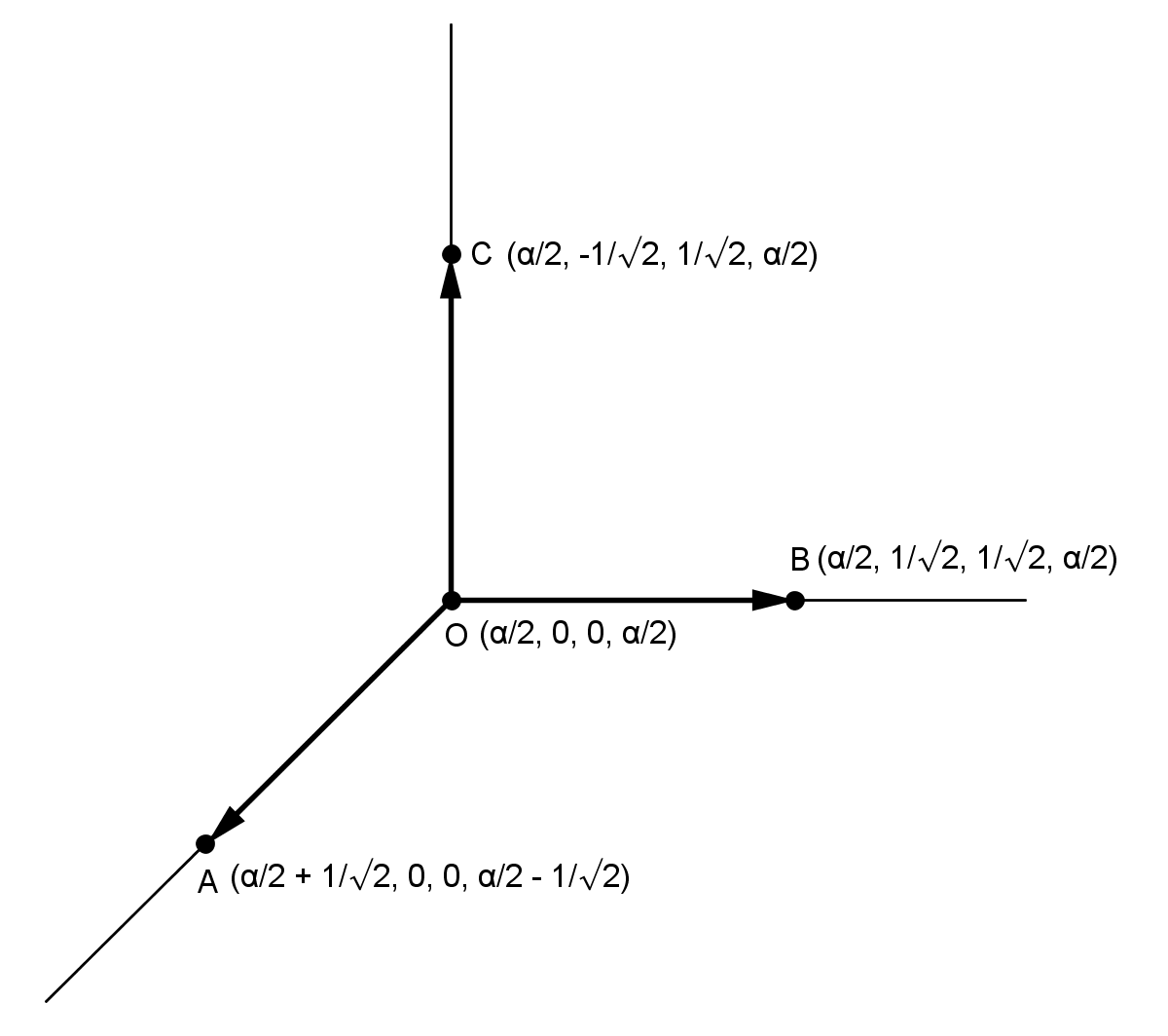}
\caption{An orthogonal Cartesian coordinate system in the hyperplane $x_1+x_4=\alpha$}\label{Fig9}
\end{center}
\end{figure}
\subsection{A coordinate system in the hyperplane $P(\alpha)$}
To describe the geometry of the set $S(\alpha, \beta)$ we introduce a rectangular cartesian coordinate system in the hyperplane $P(\alpha)$. 
Consider the following points in the hyperplane (see Figure \ref{Fig9}):
\begin{align*}
&O^\prime\left( \alpha/2, 0, 0, \alpha/2\right)\\
& A\left(\alpha/2+1/\sqrt{2}, 0, 0, \alpha/2 -1/\sqrt{2}\right)\\
& B\left( \alpha/2, 1/\sqrt{2}, 1/\sqrt{2}, \alpha/2\right)\\
& C\left( \alpha/2, -1/\sqrt{2}, 1/\sqrt{2},\alpha/2\right)
\end{align*}
These points define the following ordinary vectors in the three-dimensional space $P(\alpha)$:
\begin{align*}
\overrightarrow{O^\prime A} & = \overrightarrow{OA} - \overrightarrow{OO^\prime} \\
& = \left(\alpha/2+1/\sqrt{2}, 0, 0, \alpha/2 -1/\sqrt{2}\right) - \left( \alpha/2, 0, 0, \alpha/2\right)\\
& = \left(1/\sqrt{2}, 0, 0, -1/\sqrt{2}\right) \\
\overrightarrow{O^\prime B} & = \left( 0, 1/\sqrt{2}, 1/\sqrt{2}, 0 \right)\\
\overrightarrow{O^\prime C} & = \left( 0, - 1/\sqrt{2}, 1/\sqrt{2}, 0 \right)
\end{align*}
Using the usual inner product in $\mathbb R^4$ it can be seen that these vectors form a set of mutually perpendicular unit vectors in the three-dimensional space $P(\alpha)$. We choose $O^\prime$ as the origin and the lines $O^\prime A, O^\prime B, O^\prime C$ as the coordinate axes, the positive directions of the the coordinate axes being along the directions of the vectors $\overrightarrow{O^\prime A}, \overrightarrow{O^\prime B}, \overrightarrow{O^\prime C}$. As a tribute to R.J.T. Bell, author of a well known textbook on classical analytical three-dimensional geometry, we call this coordinate system the Bell coordinate system in $P(\alpha)$. Let $X(x_1,x_2,x_3,x_4)$ be any point in this hyperplane and let its coordinates relative to the Bell coordinate system be $(x,y,z)$. Then we have:
\begin{align*}
\overrightarrow{O^\prime X}
& = x \overrightarrow{O^\prime A} + y\overrightarrow{O^\prime B} + z \overrightarrow{O^\prime C} \\
& = x \left(1/\sqrt{2}, 0, 0, -1/\sqrt{2}\right) + 
y \left( 0, 1/\sqrt{2}, 1/\sqrt{2}, 0 \right) + 
z \left( 0, - 1/\sqrt{2}, 1/\sqrt{2}, 0 \right) \\
& = \left( x/\sqrt{2} , (y-z)/\sqrt{2},  (y+z)/\sqrt{2}, -x/\sqrt{2}\right)
\end{align*}
We also have:
\begin{align*}
\overrightarrow{O^\prime X} & = \overrightarrow{OX} - \overrightarrow{OO^\prime} \\
& = (x_1,x_2,x_3,x_4) - \left( \alpha/2, 0, 0, \alpha/2\right)\\
& = \left(x_1 - \alpha/2, x_2, x_3, x_4-\alpha/2\right)
\end{align*}

Thus we have
\begin{align*}
x_1 - \alpha/2 & = x/\sqrt{2} \\
x_2 & = (y-z)/\sqrt{2}\\
x_3 & = (y+z)/\sqrt{2}\\
x_4 - \alpha/2 & = -x/\sqrt{2}
\end{align*}
and therefore, we have
\begin{align}
x_1 & = \alpha/2 + x/\sqrt{2} \label{Eqx1}\\
x_2 & = (y-z)/\sqrt{2}\\
x_3 & = (y+z)/\sqrt{2}\\
x_4 & = \alpha/2 -x/\sqrt{2}\label{Eqx4}
\end{align}
\section{Geometry of the set $S(\alpha, \beta)$}
The following result characterises the geometry of the set $S(\alpha, \beta)$.
\begin{lemma}\qquad 
\begin{enumerate}
\item
If $\alpha^2 - 4 \beta > 0$, then $S(\alpha, \beta)$ is a hyperboloid of one sheet.
\item
If $\alpha^2 - 4 \beta = 0$, then $S(\alpha, \beta)$ is a right circular cone.
\item
If $\alpha^2 - 4 \beta < 0$, then $S(\alpha, \beta)$ is a hyperboloid of two sheets.
\end{enumerate}
\end{lemma}
\paragraph{Proof.} 
By Eq.\eqref{Eq6a} and Eq.\eqref{Eq6b}, $X=\begin{bmatrix}x_1 & x_2 \\ x_3 & x_4 \end{bmatrix} =(x_1,x_2,x_3,x_4)$ is in $S(\alpha, \beta)$ if and only if $X$ is in $P(\alpha)$ and Eq.\eqref{Eq6b} is satisfied. Now $X$ is in $P(\alpha)$ if and only if $x_1, x_2, x_3$ and $x_4$ can be expressed in the form given by Eqs.\eqref{Eqx1}-\eqref{Eqx4}, where $(x,y,z)$ are the coordinates of $X$ relative to the Bell coordinate system in $P(\alpha)$.

Thus $X$ is in $S(\alpha, \beta)$ if and only if its coordinates relative to the Bell coordinate sysytem satisfies the following equation:
\begin{equation*}
\left(\alpha/2 + x/\sqrt{2}\right) \times \left(\alpha/2 -x/\sqrt{2}\right) - \left((y-z)/\sqrt{2}\right)\times \left((y+z)/\sqrt{2}\right) = \beta
\end{equation*}
This can be simplified to 
\begin{equation}\label{Eqhy}
x^2+y^2-z^2=\frac{\alpha^2}{2}- 2 \beta.
\end{equation}
This equation represents a hyperboloid of one sheet when $\alpha^2/2 -2\beta \ne 0$, that is,  when $\alpha^2 -4\beta \ne 0$ (see \S 64 \cite{Bell}). Also Eq.\eqref{Eqhy} represents a hyperboloid of two sheets when $\alpha^2 - 4\beta<0$ (see \S 64 \cite{Bell}).  When $\alpha^2 -4\beta=0$, the equation Eq.\eqref{Eqhy} represents a right circular cone with semi-vertical angle $\frac{\pi}{4}$ (see \S 59 \cite{Bell}). $\blacksquare$
\begin{figure}
\begin{center}
\begin{tabular}{m{4cm}m{3.5cm}m{4cm}}
\centering \includegraphics[width=4.25cm, height=5cm]{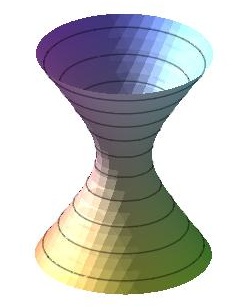} &
\centering \includegraphics[width=3cm, height=5cm]{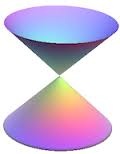} &
{\centering \includegraphics[width=3.5cm, height=5cm]{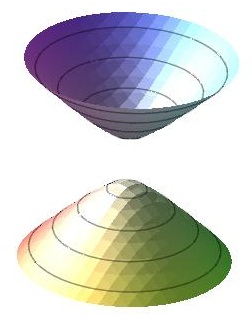}} \\
\centering (a) & \centering (b) &  \begin{center}(c)\end{center} \\ 
\centering Hyperboloid of one sheet &  \centering Right circular cone &  {\centering Hyperboloid of two sheets} \\
\end{tabular}
\caption{Images of the set $S(\alpha, \beta)$}\label{Sab}
\end{center}
\end{figure}

Figure \ref{Sab}(a) shows a sketch of a hyperboloid of one sheet, Figure \ref{Sab}(b) shows a sketch of a right circular cone and Figure \ref{Sab}(c) shows a sketch of a hyperboloid of two sheets.
\subsection{Special cases}
\begin{itemize}
\item
$S(1, 0)$ is a hyperboloid of one sheet. This has already been noted in \cite{Kris01}. This hyperboloid lies in the hyperplane $x_1+x_4=1$.
\item
$S(0,0)$ is a right circular cone with semi-vertical angle $\frac{\pi}{4}$. This fact has also been noted in \cite{Kris02}. This cone lies the hyperplane $x_1+x_4=0$. 
\item
$S(0,-1)$ is a hyperboloid of one sheet. Thus the set of square roots of the identity matrix $I_2$ defines a hyperboloid of one sheet in a four-dimensional space. This hyperboloid lies in the hyperplane $x_1+x_4=0$.
\end{itemize}
\begin{cor}
The set of all square roots of the identity matrix $I_2$ minus the two special square roots $\pm I_2$ forms a hyperboloid of one sheet in $\mathbb R^4$. The hyperboloid lies in the hyperplane $x_1+x_4=0$. 
\end{cor}
\section{The hyperboloid of square roots of $I_2$}
In this section we closely examine the hyperboloid $S(0,-1)$ formed by the
set of all square roots of $I_2$. The surface $S(0,-1)$ lies in the hyperplane $P(0)$, that is, the hyperplane represented by the equation
$$
x_1+x_4=0.
$$
To define the Bell coordinate system in this hyperplane, we choose $O^\prime(0,0,0,0)$, that is, the zero-matrix $\begin{bmatrix} 0 & 0 \\ 0 & 0 \end{bmatrix}$ as the origin. If $X(x_1,x_2,x_3,x_4)$ is any point in the plane $P(0)$ and if the Bell coordinates of $X$ are $x,y,z)$ we have
\begin{align}
x_1 & = \phantom{-1} x/\sqrt{2} \label{Eqx1b}\\
x_2 & = (y-z)/\sqrt{2}\\
x_3 & = (y+z)/\sqrt{2}\\
x_4 & = -x/\sqrt{2}\label{Eqx4a}
\end{align}
and, we also have
\begin{align}
x & = \sqrt{2}x_1 \\
y & = \sqrt{2} (x_2+x_3) \\
z & = \sqrt{2} (x_3 - x_2)\label{Eqy4}
\end{align}
The equation of $S(0,-1)$ now simplifies to
\begin{equation}
x^2+y^2-z^2 = 2.
\end{equation}
We make the following observations about this hyperboloid.
\begin{enumerate}
\item
{\bf Centre of $S(0,-1)$}

The centre of the hyperboloid $S(0,-1)$ is the origin of the Bell coordinate system in $P(0)$ which is the zero matrix in $M_2$. 
\item
{\bf Two special points  on $S(0,-1)$}

The matrix $\begin{bmatrix}1 & 0 \\ 0 & -1\end{bmatrix}$ is a square root of $I_2$. The point which corresponds to this is $(\sqrt{2},0,0)$. Similarly, 
the point which corresponds 
to $\begin{bmatrix}-1 & 0 \\ 0 & 1\end{bmatrix}$  is  $(-\sqrt{2},0,0)$. 
\item
{\bf The principal section}

The principal section of the hyperboloid $S(0,-1)$ is the intersection of the hyperboloid with the plane $z=0$. This contains all $X$ in $S(0,-1)$ which satisfies the additional condition 
$$
x_3-x_2=0\quad \text{(see Eq.\eqref{Eqy4}).}
$$
Thus the principal section of $S(0,-1)$ consists of all symmetric matrices which are square roots of $I_2$. 

Recalling the notations of Section 2, let $R=\begin{bmatrix} a & b \\ c& d \end{bmatrix}$ be an element of the principal section of $S(0,-1)$. Then we have 
\begin{align*}
a^2 + bc & =0 \\
(a+d)b & = 0 \\
(a+d)c & = 0\\
d^2 + bc & = 0 \\
b -c & =0 
\end{align*}
This system reduces to
\begin{align*}
a^2 + b^2 & = 1 \\
(a+d)b & =0 \\
d^2 - a^2 & = 0 
\end{align*}
Let $ a = \cos \phi$ and $b = \sin \phi$. If $a+d\ne 0$ then $b=0$ and $a=\pm 1$ and then $R=\pm I_2$.
If $a+d=0$ then we get
$$
R=\begin{bmatrix} \cos\phi & \sin\phi \\ \sin\phi & -\cos\phi \end{bmatrix}
$$
which is a Householder transformation. Thus, the principal section of $S(0,-1)$ represents the set of square roots of $I_2$ which are Householder transformations.
\item
{\bf The principal axis}

The principal axis is the $z$-axis specified by the equations $x=y=0$. 
These equations give the matrix $\begin{bmatrix}0 & x_2 \\ -x_2 & 0 \end{bmatrix}$ which is a skew-symmetric matrix. It follows that thw principal axis of $S(0,-1)$ is the set of skew-symmetric matrices in $M_2$.
\item
{\bf The asymptotic cone from centre}

The asymptotic cone from the centre has the equation
$$
x^2 + y^2 -z^2 = 0.
$$
This yields the equation
$$
x_1^2+4x_2x_3=0.
$$
Hence the asymptotic cone from the centre consists of all matrices  in the following set:
$$
\left\{ \begin{bmatrix}x_1 & x_2 \\ x_3 & -x_1\end{bmatrix} : x_1^2+4x_2x_3=0\right\}.
$$
\item
{\bf The generators}

One distinguishing feature of the hyperboloid of one sheet is the existence of two systems of generators on the surface. The generators are straight line lying completely on the surface. Through every point on the hyperboloid there passes one generator of each system. Let us consider the generators of the hyperboloid $S(0,-1)$. The following lemma can be used to determine the generators of this surface.
\begin{lemma}
Let $A\in S(0,-1)$ and let $U, V\in M_2$ be such that 
\begin{gather}
AU=U,\quad UA=-U,\quad U^2 =0\\
AV=-V, \quad VA=V,\quad V^2 =0.
\end{gather}
Then the two generators through $A$ are the lines specified by the sets
\begin{align*}
L_1& = \{ A+tU:t\text{ a scalar } \}\\
L_2 & = \{ A+tV: t\text{ a scalar} \}
\end{align*}
\end{lemma}
\paragraph{Proof.} We have
\begin{align*}
(A+tU)^2
& = A^2 + tAU+tUA+U^2\\
& = I_2 + tU -tU +0\\
& + I_2
\end{align*}
Hence $A+tU\in S(0,-1)$ and therefore $L_1\subseteq S(0,-1)$. Since $L_1$ is parameterized by a single parameter in the first degree, it represents a 
line. By a similar argument it can be shown that $L_2$ is also a line lying completely in $S(0,-1)$. 

Now, if possible let
$$
A+t_1U=A+t_2V.
$$
Then we have
$$
t_1U=t_2V.
$$
Pre-multiplying by $A$ we get
$$
t_1U=-t_2V.
$$
These together imply that $t_1=t_2=0$. Thus we must have $L_1\cap L_2=\emptyset$.

It follows that $L_1$ and $L_2$ are the two generators through $A$. $\blacksquare$

It is rather easy to find $U$ and $V$ satisfying the conditions of the lemma. Let $X$ be any element of $M_2$ then the following matrices have the required properties:
\begin{align*}
U& = (A+I_2)X(A-I_2)\\
V& = (A-I_2)X(A+I_2)
\end{align*}
In particular, let us find the generators through points on the principal section. A general element on the principal section is given by
$$
P=\begin{bmatrix}\cos \phi & \sin\phi \\ \sin\phi & -\cos\phi\end{bmatrix}
$$
Let us take 
$$
X=\begin{bmatrix} 1 & 0 \\ 0 & 0 \end{bmatrix}
$$
Then we get
\begin{align*}
U & =
\begin{bmatrix} -\sin^2\phi & \sin\phi(\cos\phi -1) \\ 
\sin\phi(\cos\phi+1) & \sin^2\phi \end{bmatrix} \\
V & =\begin{bmatrix} -\sin^2\phi & \sin\phi(\cos\phi +1) \\ 
\sin\phi(\cos\phi -1) & \sin^2\phi \end{bmatrix}
\end{align*}
Removing the common factor $\sin \phi$ and let it absorb with the scalar $t$, we may take 
\begin{align*}
U & = \begin{bmatrix} -\sin\phi & \cos\phi-1 \\ \cos\phi+1 & \sin\phi\end{bmatrix}\\
V & = \begin{bmatrix} -\sin\phi & \cos\phi + 1 \\ \cos\phi -1 & \sin\phi\end{bmatrix}
\end{align*}
Therefore the two generators through $P$ are the lines represented by
\begin{gather}
P+t\begin{bmatrix} -\sin\phi & \cos\phi-1 \\ \cos\phi+1 & \sin\phi\end{bmatrix}\\
P+t\begin{bmatrix} -\sin\phi & \cos\phi + 1 \\ \cos\phi -1 & \sin\phi\end{bmatrix}
\end{gather}
Incidentally, we can write $U$ as
$$
U=\begin{bmatrix}-\sin\phi & \cos\phi \\ \cos\phi & \sin\phi\end{bmatrix} + \begin{bmatrix} 0 & -1 \\ 1 & 0 \end{bmatrix}.
$$
The first term is a Householder transformation and the second term is a square root of $-I_2$. There is similar representation for $V$ also.
\end{enumerate}
\begin{figure}
\begin{center}
\includegraphics[height=5cm]{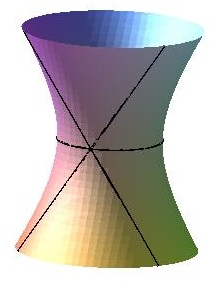}
\caption{The hyperboloid $S(0,-1)$ showing the principal section and the generators through a point on the principal section}
\end{center}
\end{figure}
\section{Matrix function approach}
\subsection{Introduction}
In this approach we begin with the concept of a function of a matrix, and then consider a definition of matrix function. We then apply this definition to the square root function of a matrix. This square root function is then specialized to compute the square roots of $I_2$.

The term ``function of a matrix'' has several meanings. Following Nigham, we look upon a function of matrix as a generalization of the concept of complex function $f(z)$ of a complex variable $z$. Accordingly, in our definition of a  matrix function, we start with some $f(z)$ and assign a meaning to $f(A)$ for any arbitrary matrix $A$. In this approach the determinant, the trace, the transpose and the adjoint are all not functions of matrices.
\subsection{Definition of function of matrix}
Let $f(z)$ be a complex function of a complex variable $z$ and let $A$ be an arbitrary square matrix of order $n$. 
\begin{itemize}
\item
Let $Z$ be a nonsingular matrix and $J$ a matrix in the Jordan canonical form such that
$$
A=ZJZ^{-1}.
$$
(Given $A$, $Z$ and $J$ exist.)
\item
Let 
$$
J=\text{diag\,}(J_1, J_2, \ldots, J_p)
$$
where
$$
J_k = J(\lambda_k) = \begin{bmatrix}
\lambda_k & 1 &  & \\
          & \lambda_k & \ddots & \\
          &           & \ddots & 1 \\
          &            &      & \lambda_k
\end{bmatrix}\in \mathbb C^{m_k\times m_k}.
$$
\item
We have
$$
m_1+m_2+\cdots+m_k=n.
$$
\item
Let $\lambda_1, \lambda_2,\ldots, \lambda_s$ be the distinct eigenvalues of $A$. Let $n_i$ be the order of the
largest Jordan block in which $\lambda_i$ appears.
\item
Assume that the values 
$$
f^{(j)}(\lambda_i)\text{ for }j=0,\ldots, n_i-1,\quad i=1,\ldots, s
$$
exist. 
\end{itemize}
Now $f(A)$ is defined as 
$$
f(A) = Z\,\,\text{diag}\, (f(J_k))Z^{-1}
$$
where we define
$$
f(J_k)=
\begin{bmatrix}
f(\lambda_k) & f^\prime(\lambda_k) & \cdots & \dfrac{f^{(m_k-1)}(\lambda_k)}{(m_k-1)!} \\
             & f(\lambda_k) & \ddots & \vdots \\
             &              & \ddots & f^\prime(\lambda_k) \\
             &              &         & f(\lambda_k)
\end{bmatrix}
$$
\subsection{Computation of square roots of $I_2$}
To compute the square roots of $I_2$ we have to consider the function
$$
f(z)=\sqrt{z}
$$
and then we have to compute the value of $f(I_2)$. For clarity of discussion, let us temporarily write $A=I_2$. The Jordan canonical form of $A$ is 
$$
J = 
\begin{bmatrix}
J_1 & 0 \\ 0 & J_2 
\end{bmatrix} 
\text{ where } J_1 = [ 1 ], \quad J_2 = [1].
$$
For any nonsingular matrix $Z$ we have
$$
A=ZJZ^{-1}.
$$
Now, by definition, we have
\begin{align*}
f(A)& = Z\, \text{diag}\,(f(J_1), f(J_2)) Z^{-1}\\
& = Z \begin{bmatrix} [\sqrt{1}] & 0 \\ 0 & [\sqrt{1}]  \end{bmatrix} Z^{-1}\\
& = Z \begin{bmatrix} \sqrt{1} & 0 \\ 0 & \sqrt{1} \end{bmatrix} Z^{-1}.
\end{align*}
The function $f(z)$ has two branches and the value $f(1)=\sqrt{1}$ can be in any one of these branches. Since $f(1)$ appears twice in $f(A)$, the two values may be chosen as lying in the same branch or in two different branches. If w choose them as lying in the same branch we get the {\em the primary values} of $f(A)$. These are given by
$$
Z\begin{bmatrix} 1 & 0 \\ 0 & 1 \end{bmatrix}Z^{-1} = I_2,\quad Z\begin{bmatrix}-1 & \phantom{-}0 \\ \phantom{-} 0 & -1 \end{bmatrix} Z^{-1}= -I_2.
$$
If we choose them as lying two different branches, we get {\em non-primary values} of $f(A)$. These are  given by
$$
Z\begin{bmatrix} \phantom{-}1 & \phantom{-}0 \\ \phantom{-}0 & - 1 \end{bmatrix}Z^{-1} ,\quad Z\begin{bmatrix}-1 & \phantom{-}0 \\ \phantom{-}0 & \phantom{-}1 \end{bmatrix} Z^{-1}.
$$
\section{Square roots of $\pm I_2$ via split-quaternions}
%
%
\subsection{Split-quaternions}
The split-quaternions, also called coquaternions, are elements of a 4-dimensional associative algebra introduced by James Cockle in 1849. They form a four dimensional real vector space equipped with a multiplication  operation. Unlike the quaternions, the algebra of split-quaternions contains zero divisors, nilpotent elements, and nontrivial idempotents. We denote the algebra of split-quaternions by $\mathsf P$.  

The set $\{1, i, j, k\}$ forms a basis for the algebra $\mathsf P$. The products of these elements are given below:
\begin{alignat*}{2}
ij & = \phantom{-}k &  & = -ji \\
jk & = -i & & = -kj\\
ki & = \phantom{-} j & & = -ik \\
i^2 & = -1 &&\\
j^2 & = +1 &&\\
k^2 & = +1 &&
\end{alignat*}
These are equivalent to the following set of equations:
$$
i^2 =-1, \quad j^2=k^2=ijk=1.
$$
The conjugate of the split-quaternion
$$
q=w+xi+yj+zk
$$
is
$$
q^*=w-xi-yj-zk
$$
and its modulus is
$$
qq^*=w^2+x^2-y^2-z^2.
$$
If $qq^*\ne 0$ then $q$ has an inverse, namely, $\dfrac{q^*}{qq^*}$.
A split-quaternion is spacelike,  lightlike or timelike according as $qq^*<0$, $qq^*=0$ or $qq^*>0$. 
\subsection{Isomorphism of $\mathsf P$ with the ring of $2\times 2$ real matrices}
Let $M_2$ be the set of real $2\times 2$ matrices. It can be easily verified that the map
$$
w+xi+yj+zk \mapsto 
\begin{bmatrix}
w+z & x+y \\ y-x & w-z
\end{bmatrix}
$$
defines a bijection from $\mathsf P$ to $M_2$ which is also a ring isomorphism.
The inverse of the map is
$$
\begin{bmatrix} a & b \\ c & d \end{bmatrix}
\mapsto
\frac{1}{2}(a+d) + \frac{1}{2}(b-c)i + \frac{1}{2}(b+c)j + \frac{1}{2}(a-d)k.
$$
If 
$$
X=\begin{bmatrix} a & b \\ c & d \end{bmatrix} \mapsto q
$$
then
$$
qq^* = ad-bc = \det(X).
$$
Also, the identity matrix $I_2$ maps to the identity $1$ in $\mathsf P$.
\subsection{Square roots of $1$ in $\mathsf P$}
Let 
$$
q=w+xi+yj+zk\in\mathbf P
$$
be such that
$$
q^2 =1.
$$
Post-multiplying by $q^*$ we have
$$
qqq^* = q^*.
$$
Clearly $qq^*\ne 0$, because otherwise we must have $q^*=0$ and so $q=0$ which is obviously a contradiction. So we must have
$$
q=\frac{q^*}{qq^*}.
$$
This yields the following equations:
\begin{align}
w\left(1 - \frac{1}{qq^*}\right)& = 0 \label{SQ1}\\
x\left(1 + \frac{1}{qq^*}\right)& = 0 \label{SQ2}\\
y\left(1 + \frac{1}{qq^*}\right)& = 0 \label{SQ3}\\
z\left(1 + \frac{1}{qq^*}\right)& = 0 \label{SQ4}
\end{align}
\subsubsection*{Case 1}

If $w\ne 0$ then $1-\frac{1}{qq^*} =0$ and so $1 + \frac{1}{qq^*}\ne 0$. Tis would then imply that $x=y=z=0$. Hence $q=w$ and $w^2=1 $ and hence $q=\pm 1 \in \mathbf P$.

\subsubsection*{Case 2}

Let $w=0$. If $1-\frac{1}{qq^*} =0$ then as above we have $x=y=z=0$ which leads to the contradiction that $q=0$. 

Now we must have 
$$
1+\frac{1}{qq^*} =0;
$$
that is, 
\begin{equation}\label{ABC}
x^2 - y^2 - z^2 = -1.
\end{equation}
Thus, the square roots of $1\in \mathsf P$, other than $\pm 1$, are split-quaternions of the form
$$
q=xi+yj+zk
$$
where $x,y,z$ satisfy Eq.\eqref{ABC}. If we identify $x,y,z$ as the coordinates of a point $(x,y,z)$ in ordinary three-dimensional euclidean space, then the set of points defined by Eq.\eqref{ABC} form a hyperboloid of one sheet. It follows that the set of square roots of $1\in \mathsf P$ can be identified with a hyperboloid of one sheet.

We can find real numbers $t$ and $\phi$ such that
\begin{align*}
x & = \sinh t \\
y & = \cosh t \sin \phi\\ 
z & = \cosh t \cos \phi
\end{align*}
Thus we have the parametrized form for the square roots of $1$:
$$
q=i \sinh t +(j \sin\phi + k\cos \phi)\cosh t
$$
\subsection{Square roots of $I_2$}
The square roots of $1\in \mathbf P$ are
$$
1, \quad -1, \quad i \sinh t +(j \sin\phi  + k\cos \phi                )\cosh t
$$
By the isomorphism between $\mathbf P$ and $M_2$, we see that the roots of $I_2$ are the following matrices:
$$
\begin{bmatrix}1 & 0 \\ 0 & 1\end{bmatrix}, \quad
\begin{bmatrix}- 1 & \phantom{-}0 \\ \phantom{-}0 & -1\end{bmatrix}, \quad
\begin{bmatrix}
 \cosh t \cos\phi  &  \cosh t \sin\phi  +\sinh t \\
\cosh t \sin \phi  - \sinh t & -  \cosh t \cos \phi  \end{bmatrix}
$$
The third form for the square root can be expressed in the following form:
$$
 \cosh t
\begin{bmatrix} \cos\phi  & \phantom{-}\sin\phi  \\ \sin\phi  & -\cos\phi  \end{bmatrix} +\sinh t \begin{bmatrix} \phantom{-}0 & 1 \\ -1 & 0 \end{bmatrix}
$$
The first term is a scalar multiple of a Holder transformation which is a square root of $I_2$ and the second term is a scalar multiple of a square root of $-I_2$. 
\subsection{Square roots of $-1\in \mathsf P$ and $-I_2$}
Using an argument similar to one used to find the square roots of $1 \in \mathsf P$, we can see  that the square roots of $-1\in \mathsf P$ are split-quaternions of the form
$$
q=xi+yj+zk
$$
where $x,y,z$ satisfy the condition
\begin{equation}\label{2S}
x^2 - y^2 -z^2 =1.
\end{equation}
Eq.\eqref{2S} represents a hyperboloid of one sheet in $\mathbb R^3$. Thus the set of square roots of $-1$ forms a hyperboloid of one sheet.

We can parameterize $x,y,z$ as follows:
\begin{align*}
x & = \sec t \\
y & = \tan t \sin \phi \\
z & = \tan t \cos \phi
\end{align*}
Then, we have
$$
q=i \sec t + j \tan t\sin\phi + k \tan t \cos\phi.
$$
Using the isomorphism between $M_2$ and $\mathsf P$ we see that  the square roots of $-I_2$ are of the form
$$
\begin{bmatrix}
\tan t \cos \phi & \tan t \sin\phi + \sec t \\
\tan t \sin\phi - \sec t & - \tan t \cos \phi
\end{bmatrix}.
$$
This can be expressed as
$$
\tan t 
\begin{bmatrix}
\cos \phi & \phantom{-}\sin\phi \\ \sin\phi & -\cos \phi 
\end{bmatrix}
+
\sec t 
\begin{bmatrix}
\phantom{-}0 & 1 \\ -1 & 0 
\end{bmatrix}
$$
The first term is a scalar multiple of a Holder transformation and the second term is a multiple of a special square root of $-I_2$.
\section{$2\times 2$ matrices with infinite number of square roots}
We have been stressing the fact that $\pm I_2$ have infinite number of square roots. Now we investigate which other $2\times 2$ matrices have infinite number of square roots. In the sequel, we shall use the following notation: For any $2\times 2$ matrix $A$, we write
$$
S(A)=\{ X \in M_2: X^2 = A\}
$$
We have to determine for which $A$'s are the sets $S(A)$ infinite. 

The following results are fairly obvious.
\begin{lemma}
For any $A\in M_2$ and any positive real number $\alpha$,
\begin{enumerate}
\item
$S(\alpha A ) = \sqrt{\alpha} S(A)$.
\item
For any nonsingular matrix $X\in M_2$, we have $S(Z^{-1}AZ)=Z^{-1}S(A)Z$.
\end{enumerate}
\end{lemma}
\paragraph{Proof.} We shall give a proof of the latter result. Let
$X\in S(Z^{-1}AZ)$ Then we have
$$
X^2 = Z^{-1}AZ.
$$
Let $Y=ZXZ^{-1}$.
\begin{align*}
Y^2 & = (ZXZ^{-1})(ZXZ^{-1})\\
& = ZX^2 Z^{-1} \\
& = Z(Z^{-1}AZ)Z^{-1}\\
& = A 
\end{align*}
Therefore, $X=Z^{-1}YZ$ with $Y\in S(A)$. This implies that 
$$
S(Z^{-1}AZ)\subseteq Z^{-1}S(A)Z.
$$
Similarly we also have
$$
Z^{-1}S(A)Z\subseteq S(Z^{-1}AZ).
$$
It follows that $S(Z^{-1}AZ)=Z^{-1}S(A)Z$. $\blacksquare$

Let $A\in M_2$ and let its eigenvalues $\lambda_1, \lambda_2$ be real.  Let $J$ be the Jordan canonical form of $A$ and let $Z$ be a non-singular   
matrix such that $A=Z^{-1}JZ$. By Lemma , 
$$
S(A)=Z^{-1}S(J)Z.
$$
Thus to find $A$'s having infinite number of square roots, it is enough to find $J$'s having infinite number of square roots. We consider the various forms of $J$ one by one.
\begin{enumerate}
\item
Let $0<\lambda_1<\lambda_2$ so that 
$$
J=
\begin{bmatrix} 
\lambda_1 & 0 \\ 0 & \lambda_2
\end{bmatrix}
$$
Clearly $S(J)$ contains only the following four elements.
$$
\begin{bmatrix}
\pm\sqrt{\lambda_1} & 0 \\
0 & \pm\sqrt{\lambda_2}
\end{bmatrix}
$$
Thus $S(A)$ is finite and contains precisely four elements.
\item
Let $0<\lambda_1=\lambda_2=\lambda$, say. 
We have to consider two sub-cases:
\begin{enumerate}
\item
Let 
$$
J =
\begin{bmatrix}
\lambda & 0 \\ 0 & \lambda 
\end{bmatrix}
=\lambda I_2
$$
Therefore by Lemma we have
$$
S(J) =\sqrt{\lambda} S(I_2).
$$
Since $S(I_2)$ is an infinite set, $S(J)$ is also an infinite set.
\item
Let
$$
J=
\begin{bmatrix}
\lambda & 1 \\ 0 & \lambda 
\end{bmatrix}
$$
Now we can find the square roots of $J$ as in Section 2. It can be shown that the following are the only two square roots of $J$ in this case:
$$
\begin{bmatrix}
\sqrt{\lambda} & 1/(2\sqrt{\lambda}) \\ 0 & \sqrt{\lambda} 
\end{bmatrix},\quad 
\begin{bmatrix}
- \sqrt{\lambda} & - 1/(2\sqrt{\lambda}) \\ 0 & - \sqrt{\lambda} \end{bmatrix}
$$
\end{enumerate}

\item
Let $0=\lambda_1 <\lambda_2 = \lambda$, say. Then we have
$$
J =
\begin{bmatrix}
\lambda & 0 \\ 0 & 0 
\end{bmatrix}
$$
In this case the only two square roots of $J$ are
$$
\begin{bmatrix}
\pm\sqrt{\lambda} & 0 \\ 0 & 0
\end{bmatrix}.
$$
\item
Let $0=\lambda_1=\lambda_2=0$. Again we have to consider two cases:
\begin{enumerate}
\item
Let
$$
J = 
\begin{bmatrix}
0 & 1 \\ 0 & 0
\end{bmatrix}.
$$
This matrix has no square roots.
\item
Let
$$
J=
\begin{bmatrix}
0 & 0 \\ 0 & 0 
\end{bmatrix}
$$
Any nilpotent matrix is a square root of $J$. Since there are infinite number of nilpotent matrices, $S(J)$ is an infinite set. 
\end{enumerate}
\end{enumerate}
\begin{lemma}
An element $A\in M_2$ with real eigenvalues has an infinite number of square roots if and only if $A=\lambda I_2$ for some $\lambda\ge 0$.
\end{lemma}
%
%
%

\end{document}